\documentclass{amsart}%
\usepackage{amsfonts}
\usepackage{amsmath}
\usepackage{amssymb}
\usepackage{graphicx}%
\newtheorem{theorem}{Theorem}
\theoremstyle{plain}

\newtheorem{corollary}[theorem]{Corollary}

\newtheorem{lemma}[theorem]{Lemma}

\newtheorem{proposition}[theorem]{Proposition}
\newtheorem{remark}[theorem]{Remark}
\numberwithin{equation}{section} \numberwithin{theorem}{section}
\begin{document}
\title[Steady states of thin film equations]{On steady states of van der Waals force driven thin film equations}
\date{\today}
\author{Huiqiang Jiang}
\address{School of Mathematics, University of Minnesota, 127 Vincent Hall, 206 Church
St. S.E., Minneapolis, MN 55455}
\email{hqjiang@math.umn.edu}
\author{Wei-Ming Ni}
\address{School of Mathematics, University of Minnesota, 127 Vincent Hall, 206 Church
St. S.E., Minneapolis, MN 55455}
\email{ni@math.umn.edu}

\begin{abstract}
Let $\Omega\subset\mathbb{R}^{N}$, $N\geq2$ be a bounded smooth domain and
$\alpha>1$. We are interested in the singular elliptic equation%
\[
\triangle h=\frac{1}{\alpha}h^{-\alpha}-p\quad\text{in }\Omega
\]
with Neumann boundary conditions. In this paper, we gave a complete
description of all continuous radially symmetric solutions. In particular, we
constructed nontrivial smooth solutions as well as rupture solutions. Here a
continuous solution is said to be a rupture solution if its zero set is
nonempty. When $N=2$ and $\alpha=3$, the equation has been used to model
steady states of van der Waals force driven thin films of viscous fluids. We
also considered the physical problem when total volume of the fluid is prescribed.

\end{abstract}
\maketitle

\section{Introduction}

The equation%
\begin{equation}
h_{t}=\triangledown\cdot\left(  h^{3}\triangledown p\right)
\label{Equation VdW dynamics}%
\end{equation}
has been used to model the dynamics of van der Waals force driven thin films
of viscous fluids\cite{ntofr}\cite{MR1707809}\cite{dotdtfr}\cite{MR1707811}.
Here $h$ is the thickness of the thin film and the pressure
\begin{equation}
p=\frac{1}{3}h^{-3}-\triangle h, \label{Equation VdW pressure}%
\end{equation}
is a sum of contributions from disjoining pressure due to attractive van der
Waals force and a linearized curvature term corresponding to surface tension
effects. Hence, $\left(  \ref{Equation VdW dynamics}\right)  $ becomes%
\begin{equation}
h_{t}=-\triangledown\cdot\left(  h^{-1}\triangledown h\right)  -\triangledown
\cdot\left(  h^{3}\triangledown\triangle h\right)  ,
\label{Equation dynamics vdw}%
\end{equation}
which is a special case of the generalized thin film equation%
\begin{equation}
h_{t}=-\triangledown\cdot\left(  h^{m}\triangledown h\right)  -\triangledown
\cdot\left(  h^{n}\triangledown\triangle h\right)  \label{Equation dynamics}%
\end{equation}
where the exponents $m,n$ represent the powers in the destabilizing
second-order and the stabilizing fourth-order diffusive terms, respectively.
This class of equations occurs in connection with many physical models
involving fluid interfaces\cite{MR1642807}\cite{lseotlf}. For example, when
$n=1$ and $m=1$, it describes a gravity driven Hele-Shaw cell\cite{MR1396471}%
\cite{MR1277195}\cite{MR1377907}\cite{MR1377908}\cite{MR1864368} ; for $n=m=3$
it describes fluid droplets hanging from a ceiling\cite{tsohd}; and for $n=0$
and $m=1$, it is a modified Kuramoto-Sivashinsky equation which describes
solidification of a hyper-cooled melt\cite{MR1345326}\cite{MR1836532}. Over
the past two decades, these models have also been the focus of rigorous and
extensive mathematical analysis\cite{MR1471465}\cite{MR1455063}%
\cite{MR1371925}\cite{MR1627165}\cite{MR1751951}\cite{MR1616558}
\cite{2004Grun}\cite{MR1778120}\cite{MR1844589}\cite{MR1900328}%
\cite{MR1938391}\cite{2005Slepcev_Pugh}.

As in the van der Waals force case, when $n-m\neq1$, letting%
\begin{equation}
p=-\frac{1}{m-n+1}h^{m-n+1}-\triangle h, \label{Equation pressure mn}%
\end{equation}
we can rewrite $\left(  \ref{Equation dynamics}\right)  $ as
\[
h_{t}=\triangledown\cdot\left(  h^{n}\triangledown p\right)  .
\]

Now we consider viscous fluids in a cylindrical container whose bottom is
represented by $\Omega$, a bounded smooth domain in $\mathbb{R}^{2}$. Since
there is no flux across the boundary, we have the Neumann boundary condition%
\begin{equation}
\frac{\partial p}{\partial\nu}=0\text{ on }\partial\Omega.
\label{Equation no flux}%
\end{equation}
We also ignore the wetting or nonwetting effect, and assume that the fluid
surface is perpendicular to the boundary of the container, i.e.,
\begin{equation}
\frac{\partial h}{\partial\nu}=0\text{ on }\partial\Omega.
\label{Equation nonwetting}%
\end{equation}

Whenever $m-n\neq-1$ or $-2$, we can associate $\left(
\ref{Equation dynamics}\right)  $ with energy%
\[
E\left(  h\right)  =\int_{\Omega}\left(  \frac{1}{2}\left\vert \triangledown
h\right\vert ^{2}-\frac{1}{\left(  m-n+1\right)  \left(  m-n+2\right)
}h^{m-n+2}\right)  ,
\]
and formally, using $\left(  \ref{Equation no flux}\right)  ,\left(
\ref{Equation nonwetting}\right)  $, we have%
\begin{align*}
&  \frac{d}{dt}E\left(  h\right)  =\int_{\Omega}\left(  -\triangle
hh_{t}-\frac{1}{m-n+1}h^{m-n+1}h_{t}\right) \\
=  &  \int_{\Omega}p\triangledown\cdot\left(  h^{n}\triangledown p\right)
=-\int_{\Omega}h^{n}\left\vert \triangledown p\right\vert ^{2}.
\end{align*}
Hence, for a thin film fluid at rest, $p$ has to be a constant, and $h$
satisfies $\left(  \ref{Equation pressure mn}\right)  $.

Therefore, letting $\alpha=-\left(  m-n+1\right)  $, we are led to the
elliptic problem%
\begin{equation}
\left\{
\begin{array}
[c]{lll}%
\triangle h=\frac{1}{\alpha}h^{-\alpha}-p & \text{in} & \Omega,\\
\frac{\partial h}{\partial\nu}=0 & \text{on} & \partial\Omega,
\end{array}
\right.  \label{Equation prescribe pressure}%
\end{equation}
where $\Omega$ is a bounded smooth domain in $\mathbb{R}^{N}$, $N\geq1$ and
$p$ is a constant.

When $N=1$, this equation has been studied by R. Laugesen and M. Pugh in
\cite{MR1844589} where they produced positive, smooth steady states for all
$\alpha$ and touchdown steady states for $\alpha<1$. In
\cite{2001Bertozzi_Grun_Witelski}, A. L. Bertozzi, G. Gr\"{u}n and T. P.
Witelski considered $\left(  \ref{Equation prescribe pressure}\right)  $ with
additional Born repulsion term which leads to the elliptic equation%
\begin{equation}
\left\{
\begin{array}
[c]{lll}%
\triangle h=\frac{1}{\alpha}h^{-\alpha}\left(  1-\left(  \frac{\varepsilon}%
{h}\right)  ^{\beta}\right)  -p & \text{in} & \Omega,\\
\frac{\partial h}{\partial\nu}=0 & \text{on} & \partial\Omega,
\end{array}
\right.  \label{Equation Born repulsion}%
\end{equation}
where $\beta$ is a positive constant. When $\varepsilon>0$, the associated
energy to $\left(  \ref{Equation Born repulsion}\right)  $ is bounded from
below which makes a variational approach possible and enables them to show the
existence of an energy minimizer in any dimensions. It seems difficult to
extend this approach to the limiting case $\varepsilon=0$.

The goal of this paper is to understand radial solutions of $\left(
\ref{Equation prescribe pressure}\right)  $ when $N\geq2$ and $\alpha>1$ which
we will assume throughout this paper. In particular, when $N=2$ and $\alpha
=3$, we come to the van der Waals force driven thin films in the physically
realistic dimension. When $\alpha>1$, except the limited discussions in
\cite{jiang_lin}, there seems no established elliptic theory for $\left(
\ref{Equation prescribe pressure}\right)  $, and hence it is the
mathematically more interesting case. We remark that energy method can be
applied to yield nontrivial solutions to $\left(
\ref{Equation prescribe pressure}\right)  $ when $\alpha<1$. On the other
hand, the behavior of radial solutions is also quite different when $\alpha
>1$. For example, we will show that for $N\geq2$, the radial solutions will
never vanish away from the origin which contrasts with the $\alpha<1$ case
where touchdown steady states can be shown to exist in any dimensions.

Due to the singular nature of $\left(  \ref{Equation prescribe pressure}%
\right)  $, we need to be careful in discussing "solutions" to $\left(
\ref{Equation prescribe pressure}\right)  $. We say $h$ is a
\textit{continuous solution} of $\left(  \ref{Equation prescribe pressure}%
\right)  $ in $\Omega$, if $h\not \equiv 0$ and is a nonnegative continuous
function in $\overline{\Omega}$ satisfying the equation in $\left(
\ref{Equation prescribe pressure}\right)  $ in the open set $\left\{
x\in\Omega:h\left(  x\right)  >0\right\}  $. The \textit{rupture set} of $h$,%
\[
\Sigma=\left\{  x\in\Omega:h\left(  x\right)  =0\right\}  ,
\]
corresponds to "dry spots" in the thin film, which is of great significance in
the coatings industry where nonuniformities are very undesirable. Standard
elliptic theory implies that $h$ is smooth and hence a classical solution of
$\left(  \ref{Equation prescribe pressure}\right)  $ in $\Omega\backslash
\Sigma$. An interesting Hausdorff dimension estimate of $\Sigma$ can be found
in \cite{jiang_lin} where it is shown that any finite energy solution
satisfies $\mathcal{H}^{\mu}\left(  \Sigma\right)  =0$ where $\mu=N-2+\frac
{4}{\alpha+1}$. For van der Waals force driven thin film, we have $N=2$ and
$\alpha=3$, hence $\mathcal{H}^{1}\left(  \Sigma\right)  =0$, i.e., the thin
film with finite energy can't have one dimensional rupture set.

For any $p>0$, let%
\begin{equation}
\xi=\left(  \alpha p\right)  ^{-\frac{1}{\alpha}} \label{Equation xi}%
\end{equation}
then $h\equiv\xi$ is always a solution to $\left(
\ref{Theorem prescribe volume}\right)  $. The natural question is whether it
is the only solution. For the radially symmetric case, after a simple scaling,
the uniqueness theorem in \cite{delpino} implies:

\begin{proposition}
\label{Proposition uniqueness}Let $N\geq2$ and $\alpha>1$. For any given
$R>0$, there exists a constant $p_{0}$, such that for any $p\leq p_{0}$,
$h\equiv\xi$ is the only radial solution of $\left(
\ref{Equation prescribe pressure}\right)  $ in $B_{R}\left(  0\right)  $.
\end{proposition}

When $p$ is large, nontrivial solutions do exist. In fact, we have

\begin{theorem}
\label{Theorem pk}Let $N\geq2$ and $\alpha>1$. For any given $R>0$, there
exists a nondecreasing sequence of $\left\{  p_{k}\right\}  $, such that for
any $p>p_{k}$, $\left(  \ref{Equation prescribe pressure}\right)  $ has at
least $k$ nontrivial smooth radial solutions in $B_{R}\left(  0\right)  $.
\end{theorem}

Theorem \ref{Theorem pk} is an application of Theorem \ref{Theorem smooth}
below which gives a complete description of all nontrivial smooth radial solutions.

\begin{theorem}
\label{Theorem smooth}Let $N\geq2$ and $\alpha>1$. For any given $p>0$, and
for any $\eta>0$, $\eta\neq\xi$, there exists an increasing sequence $\left\{
r_{k}^{p,\eta}\right\}  _{k=1}^{\infty}$ with%
\[
r_{1}^{p,\eta}\geq\max\left\{  \sqrt{\frac{2N\alpha\left(  \eta-\xi\right)
}{\xi^{-\alpha}-\eta^{-\alpha}}},\left(  \frac{p_{0}}{p}\right)
^{\frac{1+\alpha}{2\alpha}}\right\}
\]
where $p_{0}$ is the constant in Proposition \ref{Proposition uniqueness}, and%
\begin{equation}
\lim_{k\rightarrow\infty}\left(  r_{k+1}^{p,\eta}-r_{k}^{p,\eta}\right)
=\pi\left(  \alpha p\right)  ^{-\frac{1+\alpha}{2\alpha}},
\end{equation}
such that for each $r_{k}^{p,\eta}$, there exists a unique smooth radial
solution of $\left(  \ref{Equation prescribe pressure}\right)  $ in
$B_{r_{k}^{p,\eta}}\left(  0\right)  $ satisfying $h\left(  0\right)  =\eta$.
\end{theorem}

We say a continuous solution to $\left(  \ref{Equation prescribe pressure}%
\right)  $ is a \textit{rupture solution} if $\Sigma$ is not empty. It will be
shown that for radial solutions, rupture can only occur at the origin. (See
Corollary \ref{Corollary global} below.) Our main result in this paper is as follows.

\begin{theorem}
\label{Theorem rupture}Let $N\geq2$ and $\alpha>1$. For any given $p>0$, there
exists an increasing sequence $\left\{  r_{k}^{p,0}\right\}  _{k=1}^{\infty}$
with%
\begin{equation}
r_{1}^{p,0}\geq\left(  \frac{p_{0}}{p}\right)  ^{\frac{1+\alpha}{2\alpha}},
\end{equation}
and%
\begin{equation}
\lim_{k\rightarrow\infty}\left(  r_{k+1}^{p,0}-r_{k}^{p,0}\right)  =\pi\left(
\alpha p\right)  ^{-\frac{1+\alpha}{2\alpha}},
\end{equation}
where $p_{0}$ is the constant in Proposition \ref{Proposition uniqueness},
such that for each $r_{k}^{p,0}$, there exists a unique radial rupture
solution to $\left(  \ref{Equation prescribe pressure}\right)  $ in
$B_{r_{k}^{p,0}}\left(  0\right)  $. Furthermore, if $R\neq r_{k}^{p,0}$ for
any $k$, then there is no radial rupture solution to $\left(
\ref{Equation prescribe pressure}\right)  $ in $B_{R}\left(  0\right)  $.
\end{theorem}

\begin{remark}
The rupture solutions constructed in Theorem \ref{Theorem rupture} above are
weak solutions to $\left(  \ref{Equation prescribe pressure}\right)  $ in the
distributional sense. (See Remark \ref{Remark weak} below.) We also remark
that when $N=1$ and $\alpha>1$, there is no radial rupture solutions
\cite{MR1844589}.
\end{remark}

For any given $p>0$, there exists a nontrivial radial solutions to $\left(
\ref{Equation prescribe pressure}\right)  $ in $B_{R}\left(  0\right)  $ if
and only if $r_{k}^{p,\eta}=R$ holds for some $\eta\geq0$, $\eta\neq\xi$ and
for some integer $k\geq1$.

In physical experiments, usually the total volume of the fluid is known, i.e.,
the average film thickness%
\[
\bar{h}=\frac{1}{\left\vert \Omega\right\vert }\int_{\Omega}h\left(  x\right)
dx
\]
is given while the pressure $p$ is an unknown constant. Hence, given $\bar
{h}>0$, we need to find function $h$ and constant $p$, such that
\begin{equation}
\left\{
\begin{array}
[c]{lll}%
\triangle h=\frac{1}{\alpha}h^{-\alpha}-p & \text{in} & \Omega,\\
\frac{1}{\left\vert \Omega\right\vert }\int_{\Omega}h\left(  x\right)
dx=\bar{h} &  & \\
\frac{\partial h}{\partial\nu}=0 & \text{on} & \partial\Omega.
\end{array}
\right.  \label{Equation prescribe volume}%
\end{equation}
When $\Omega=B_{1}\left(  0\right)  $, all radial solutions of $\left(
\ref{Equation prescribe volume}\right)  $ can be obtained by scaling from
solutions in Theorems \ref{Theorem smooth} and \ref{Theorem rupture}. We will
discuss such scaling in Section \ref{Section scaling}. In particular, we will show

\begin{theorem}
\label{Theorem prescribe volume}Let $N\geq2$, $\alpha>1$ and $\Omega
=B_{1}\left(  0\right)  \subset\mathbb{R}^{N}$. There exists a sequence of
thickness $\bar{h}_{1},\bar{h}_{2},\cdots$ satisfying%
\[
\lim_{k\rightarrow\infty}\sqrt{k\pi}\bar{h}_{k}=1,
\]
such that for any $k$, $\left(  \ref{Equation prescribe volume}\right)  $ with
$\bar{h}=\bar{h}_{k}$ has a radial rupture solution which, viewed as a
function in $r$, has exactly $k-1$ critical points in $\left(  0,1\right)  $.
Furthermore, if $\bar{h}\neq\bar{h}_{k}$ for any $k$, then $\left(
\ref{Equation prescribe volume}\right)  $ has no radial rupture solution.
\end{theorem}

When $\Omega=B_{1}\left(  0\right)  $, Proposition
\ref{Proposition uniqueness} implies that nontrivial solutions to $\left(
\ref{Equation prescribe volume}\right)  $ must satisfy $p>p_{0}$. Since%
\[
p=\frac{1}{\left\vert \Omega\right\vert }\int_{\Omega}\frac{1}{\alpha
}h^{-\alpha},
\]
we may ask the existence of a critical average film thickness $\bar{h}_{0}$ so
that there is no nontrivial solutions to $\left(
\ref{Equation prescribe volume}\right)  $ whenever $h\geq\bar{h}_{0}$.
Numerical analysis suggests that $\bar{h}_{0}$ does not exist. However, we are
unable to provide an analytical proof. Such a proof could be possible if we
have better understanding of $\bar{h}\left(  p,\eta,k\right)  $ which is
defined in section 5. Moreover, the detailed property of $\bar{h}\left(
p,\eta,k\right)  $ could also provide us a statement similar to Theorem
\ref{Theorem prescribe volume} for the smooth radial solutions.

The paper is organized in the following way: In Section 2, we show that any
radial solution can be extended to a global solution which is oscillating
around $\xi$. In Section 3 and Section 4, we discuss smooth radial solutions
and rupture solutions respectively and Theorems \ref{Theorem smooth} and
\ref{Theorem rupture} are proved. In Section 5, we use scaling argument to
prove Theorems \ref{Theorem pk} and \ref{Theorem prescribe volume}.

\section{Preliminaries}

Recall that given $\alpha>1$ and $p>0$, $h\in C^{0}\left(  B_{R}\left(
0\right)  \right)  $ is said to be a continuous solution of%
\begin{equation}
\triangle h=\frac{1}{\alpha}h^{-\alpha}-p \label{Equation main}%
\end{equation}
in $B_{R}\left(  0\right)  $ if $h\geq0$ and it satisfies $\left(
\ref{Equation main}\right)  $ in the open set $\left\{  x\in B_{R}\left(
0\right)  :h\left(  x\right)  >0\right\}  $.

Let $h$ be a radially symmetric solution to $\left(  \ref{Equation main}%
\right)  $, we can view $h$ as a continuous function defined on $\left[
0,R\right)  $ satisfying%
\begin{equation}
h^{\prime\prime}+\frac{N-1}{r}h^{\prime}+f\left(  h\right)  =0
\label{Equation radially symmetric}%
\end{equation}
in the set $S^{+}=\left\{  r\in\left(  0,R\right)  :h\left(  r\right)
>0\right\}  $. Here%
\[
f\left(  h\right)  =-\frac{1}{\alpha}h^{-\alpha}+p
\]
is monotone increasing and its antiderivative%
\[
F\left(  h\right)  =\frac{1}{\alpha\left(  \alpha-1\right)  }h^{1-\alpha}+ph
\]
is convex in $\left(  0,\infty\right)  $. Let
\[
\xi=\left(  \alpha p\right)  ^{-\frac{1}{\alpha}},
\]
then $F^{\prime}\left(  \xi\right)  =f\left(  \xi\right)  =0$, and $F$
achieves its absolute minimum at $\xi$. Furthermore, $F\left(  h\right)
\rightarrow\infty$ as $h\rightarrow0^{+}$ or $h\rightarrow\infty$.

For each $r\in S^{+}$, letting%
\begin{align*}
e_{1}\left(  r\right)   &  =\frac{1}{2}\left(  h^{\prime}\left(  r\right)
\right)  ^{2}+F\left(  h\left(  r\right)  \right)  ,\\
e_{2}\left(  r\right)   &  =\frac{1}{2}\left(  r^{N-1}h^{\prime}\left(
r\right)  \right)  ^{2}+r^{2\left(  N-1\right)  }F\left(  h\left(  r\right)
\right)  =r^{2\left(  N-1\right)  }e_{1}\left(  r\right)  ,
\end{align*}
we have%
\begin{equation}
\frac{d}{dr}\left[  e_{1}\left(  r\right)  \right]  =-\frac{N-1}{r}\left(
h^{\prime}\left(  r\right)  \right)  ^{2}\leq0, \label{Equation e1}%
\end{equation}
and%
\begin{equation}
\frac{d}{dr}\left[  e_{2}\left(  r\right)  \right]  =2\left(  N-1\right)
F\left(  h\left(  r\right)  \right)  r^{2N-3}\geq0.
\end{equation}

The monotonicity of $e_{1}$ and $e_{2}$ will be used to obtain a priori bounds.

\begin{lemma}
\label{Lemma no rupture}Let $0<r_{1}<r_{2}$ be such that $\left(  r_{1}%
,r_{2}\right)  \subset S^{+}$. Given $\bar{r}\in\left(  r_{1},r_{2}\right)  $,
we have, for any $r\in\left(  r_{1},r_{2}\right)  $,%
\begin{equation}
e_{1}\left(  r\right)  \leq\left(  \frac{\bar{r}}{r}\right)  ^{2\left(
N-1\right)  }e_{1}\left(  \bar{r}\right)  . \label{Equation bounds on e_1}%
\end{equation}
Furthermore,
\begin{equation}
c_{1}\leq h\left(  r\right)  \leq c_{2} \label{Equation bounds on h}%
\end{equation}
where $c_{1},c_{2}$ are two positive constants depending on $\alpha
,p,N,r_{1},\bar{r},h\left(  \bar{r}\right)  $ and $h^{\prime}\left(  \bar
{r}\right)  $, and are independent of $r_{2}$.
\end{lemma}

\begin{proof}
Since $e_{2}\left(  r\right)  $ is monotone increasing, we have, for any
$r\in\left(  r_{1},\bar{r}\right]  $,%
\begin{equation}
e_{1}\left(  r\right)  =r^{2\left(  1-N\right)  }e_{2}\left(  r\right)  \leq
r^{2\left(  1-N\right)  }e_{2}\left(  \bar{r}\right)  =\left(  \frac{\bar{r}%
}{r}\right)  ^{2\left(  N-1\right)  }e_{1}\left(  \bar{r}\right)  .
\label{Equation bounds on e1 1}%
\end{equation}
On the other hand, since $e_{1}\left(  r\right)  $ is monotone decreasing, we
have, for any $r\in\left[  \bar{r},r_{2}\right)  $,%
\begin{equation}
e_{1}\left(  r\right)  \leq e_{1}\left(  \bar{r}\right)  .
\label{Equation bounds on e1 2}%
\end{equation}
Combining $\left(  \ref{Equation bounds on e1 1}\right)  $ and $\left(
\ref{Equation bounds on e1 2}\right)  $, we obtain $\left(
\ref{Equation bounds on e_1}\right)  $. Now for any $r\in\left(  r_{1}%
,r_{2}\right)  $,
\[
F\left(  h\left(  r\right)  \right)  \leq e_{1}\left(  r\right)  \leq\left(
\frac{\bar{r}}{r_{1}}\right)  ^{2\left(  N-1\right)  }e_{1}\left(  \bar
{r}\right)  ,
\]
so $\left(  \ref{Equation bounds on h}\right)  $ follows from the fact that
$F\left(  h\right)  \rightarrow\infty$ as $h\rightarrow0^{+}$ or
$h\rightarrow\infty$.
\end{proof}

\begin{corollary}
\label{Corollary global}$h$ can not have rupture away from the origin, i.e.,
$S^{+}=\left(  0,R\right)  $. Furthermore, $h$ can be uniquely extended to a
positive smooth solution of $\left(  \ref{Equation radially symmetric}\right)
$ in $\left(  0,\infty\right)  $.
\end{corollary}

\begin{proof}
Since $S^{+}$ is open, it is a union of open intervals of the form $\left(
r_{1},r_{2}\right)  $ with $r_{1},r_{2}\notin S^{+}$. Given any such interval,
if $r_{1}>0$, Lemma \ref{Lemma no rupture} implies
\[
\liminf_{r\rightarrow r_{1}^{+}}h\left(  r\right)  >0,
\]
and since $h$ is continuous, we conclude $h\left(  r_{1}\right)  >0$, which
contradicts the assumption $r_{1}\notin S^{+}$. Similarly, we can get a
contradiction if $r_{2}<R$. Hence, $S^{+}=\left(  0,R\right)  $. Extending $h$
to a maximal interval of existence $\left(  0,R^{\ast}\right)  $. If $R^{\ast
}<\infty$, applying Lemma \ref{Lemma no rupture} again, we have for some
positive constants $c_{1},c_{2}$%
\[
c_{1}\leq h\left(  r\right)  \leq c_{2}\text{ for any }r\in\left(
R/2,R^{\ast}\right)  ,
\]
so the solution can be extended beyond $R^{\ast}$. Hence, $R^{\ast}=\infty$.
\end{proof}

Now, redefining $S^{+}=\left\{  r>0:h\left(  r\right)  >0\right\}  $, we
observe that $S^{+}=\left(  0,\infty\right)  $ and Lemma
\ref{Lemma no rupture} still holds. In particular, $\left(
\ref{Equation bounds on h}\right)  $ holds for all $r_{1}<r<\infty$. In the
remaining part of this section, we shall show that $h$ oscillates around $\xi$
near $r=\infty$.

We will need Sturm's Separation Theorem.

\begin{lemma}
\label{Lemma Sturm}\cite{ode} Let $q\left(  t\right)  $ be a real-valued
continuous function such that%
\[
0<m\leq q\left(  t\right)  \leq M.
\]
Given $t_{2}>t_{1}>0$, if $u=u\left(  t\right)  $ is a nontrivial solution of%
\[
u^{\prime\prime}+q\left(  t\right)  u=0
\]
satisfying $u\left(  t\right)  >0$ on $\left(  t_{1},t_{2}\right)  $, then%
\[
t_{2}-t_{1}\leq\frac{\pi}{\sqrt{m}}.
\]
And if in addition $u\left(  t_{1}\right)  =u\left(  t_{2}\right)  =0$, then%
\[
t_{2}-t_{1}\geq\frac{\pi}{\sqrt{M}}.
\]

\end{lemma}

\begin{lemma}
\label{Lemma existence of critical points}For any $r_{0}>0$, there exists
$r_{1}>r_{0}$ such that $h^{\prime}\left(  r_{1}\right)  =0$.
\end{lemma}

\begin{proof}
Suppose this is false, then we have either $h^{\prime}\left(  r\right)  >0$
for all $r\in\left(  r_{0},\infty\right)  $ or $h^{\prime}\left(  r\right)
<0$ for all $r\in\left(  r_{0},\infty\right)  $. Hence, $h$ is strictly
monotone increasing or decreasing on $\left(  r_{0},\infty\right)  $. From
Lemma $\ref{Lemma no rupture}$ and the observation above, it follows that $h$
is also bounded at $\infty$. So we can assume%
\[
\lim_{r\rightarrow\infty}h\left(  r\right)  =\zeta
\]
for some $\zeta>0$. For any $r>r_{0}$, integrating $\left(
\ref{Equation radially symmetric}\right)  $ from $r_{0}$ to $r$, we obtain%
\begin{equation}
h^{\prime}\left(  r\right)  -h^{\prime}\left(  r_{0}\right)  +\int_{r_{0}}%
^{r}\frac{N-1}{s}h^{\prime}\left(  s\right)  ds+\int_{r_{0}}^{r}f\left(
h\left(  s\right)  \right)  ds=0. \label{Equation difference of
derivatives}%
\end{equation}
Since%
\[
\frac{1}{2}\left(  h^{\prime}\left(  r\right)  \right)  ^{2}\leq e_{1}\left(
r\right)  \leq e_{1}\left(  r_{0}\right)  ,
\]
$h^{\prime}\left(  r\right)  $ is bounded in $\left[  r_{0},\infty\right)  $.
Now, as $h^{\prime}\left(  r\right)  $ does not change sign,
\[
\left\vert \int_{r_{0}}^{r}\frac{N-1}{s}h^{\prime}\left(  s\right)
ds\right\vert \leq\frac{N-1}{r_{0}}\left\vert \int_{r_{0}}^{r}h^{\prime
}\left(  s\right)  ds\right\vert =\frac{N-1}{r_{0}}\left\vert h\left(
r\right)  -h\left(  r_{0}\right)  \right\vert
\]
which is also bounded in $\left[  r_{0},\infty\right)  $. Hence identity
$\left(  \ref{Equation difference of derivatives}\right)  $ implies that%
\[
\int_{r_{0}}^{r}f\left(  h\left(  s\right)  \right)  ds
\]
is bounded in $\left[  r_{0},\infty\right)  $. Since $\lim_{r\rightarrow
\infty}h\left(  r\right)  =\zeta$, we have%
\[
\lim_{r\rightarrow\infty}f\left(  h\left(  r\right)  \right)  =f\left(
\zeta\right)
\]
which must be $0$. Thus $\zeta=\xi$ and%
\[
\lim_{r\rightarrow\infty}h\left(  r\right)  =\xi\text{.}%
\]
Now let%
\[
v\left(  r\right)  =r^{\frac{N-1}{2}}\left(  h\left(  r\right)  -\xi\right)
\text{,}%
\]
then%
\begin{equation}
v_{rr}+\frac{\left(  N-1\right)  \left(  3-N\right)  }{4r^{2}}v+r^{\frac
{N-1}{2}}f\left(  h\right)  =0. \label{Equation scaled difference of h and Xi}%
\end{equation}
Since $f\left(  \xi\right)  =0$, we can rewrite $\left(
\ref{Equation scaled difference of h and Xi}\right)  $ as%
\begin{equation}
v_{rr}+B\left(  r\right)  v=0,
\label{Equation scaled difference of h and Xi 2}%
\end{equation}
where%
\[
B\left(  r\right)  =\frac{f\left(  h\right)  -f\left(  \xi\right)  }{h-\xi
}+\frac{\left(  N-1\right)  \left(  3-N\right)  }{4r^{2}}.
\]
Now $\lim_{r\rightarrow\infty}h\left(  r\right)  =\xi$ implies%
\[
\lim_{r\rightarrow\infty}B\left(  r\right)  =f^{\prime}\left(  \xi\right)
>0.
\]
By Lemma \ref{Lemma Sturm},%
\[
v\left(  r\right)  =r^{\frac{N-1}{2}}\left(  h\left(  r\right)  -\xi\right)
\]
will be oscillating around $0$ as $r\rightarrow\infty$, which contradicts the
assumption that $h\left(  r\right)  \rightarrow\xi$ in a strictly monotonic manner.
\end{proof}

Next, we have

\begin{lemma}
\label{Lemma global solution is oscillating}Let $h^{\prime}\left(
r_{0}\right)  =0$ for some $r_{0}\geq0$.\newline(i). If $h\left(
r_{0}\right)  =\xi$, then $h\left(  r\right)  \equiv\xi$.\newline(ii). If
$h\left(  r_{0}\right)  >\xi$, then there exists $r_{1}>r_{0}$, such that
$h^{\prime}\left(  r\right)  <0$ on $\left(  r_{0},r_{1}\right)  $,
$h^{\prime}\left(  r_{1}\right)  =0$, $h\left(  r_{1}\right)  <\xi$ and
$F\left(  h\left(  r_{1}\right)  \right)  <F\left(  h\left(  r_{0}\right)
\right)  $.\newline(iii). If $0<h\left(  r_{0}\right)  <\xi$, then there
exists $r_{1}>r_{0}$, such that $h^{\prime}\left(  r\right)  >0$ on $\left(
r_{0},r_{1}\right)  $, $h^{\prime}\left(  r_{1}\right)  =0$, $h\left(
r_{1}\right)  >\xi$ and $F\left(  h\left(  r_{1}\right)  \right)  <F\left(
h\left(  r_{0}\right)  \right)  $.
\end{lemma}

\begin{proof}
(i). This is the standard ODE uniqueness result.

(ii). Since $h\left(  r_{0}\right)  >\xi$, we have $f\left(  h\left(
r_{0}\right)  \right)  >0$. Now%
\[
\left(  r^{N-1}h^{\prime}\right)  ^{\prime}=-r^{N-1}f\left(  h\right)  ,
\]
implies that $r^{N-1}h^{\prime}$ is strictly monotone decreasing in $\left(
r_{0},r_{0}+\delta\right)  $ for some $\delta>0$. Hence we have $h^{\prime
}\left(  r\right)  <0$ on $\left(  r_{0},r_{0}+\delta\right)  $. Applying
Lemma \ref{Lemma existence of critical points}, there exists $r_{1}>r_{0}$,
such that $h^{\prime}\left(  r_{1}\right)  =0$, and we also have $h^{\prime
}\left(  r\right)  <0$ on $\left(  r_{0},r_{1}\right)  $ if we choose the
smallest such $r_{1}$. If $h\left(  r_{1}\right)  >\xi$, we would have
$r^{N-1}h^{\prime}$ is strictly decreasing near $r_{1}$, hence $h^{\prime
}\left(  r_{1}\right)  <0$, which gives a contradiction. And if $h\left(
r_{1}\right)  =\xi$, then $h\equiv\xi$, which contradicts the hypothesis
$h\left(  r_{0}\right)  >\xi$. Hence we have $h\left(  r_{1}\right)  <\xi$.
Finally, $F\left(  h\left(  r_{1}\right)  \right)  <F\left(  h\left(
r_{0}\right)  \right)  $ follows from $\left(  \ref{Equation e1}\right)  $.

(iii). Similar to the proof of part (ii).
\end{proof}

Let $h$ be a nontrivial global solution of $\left(
\ref{Equation radially symmetric}\right)  $, starting with $r_{1}>0$ such that
$h^{\prime}\left(  r_{1}\right)  =0$. The existence of $r_{1}$ is guaranteed
by Lemma \ref{Lemma existence of critical points}. Without loss of generality,
we assume $h\left(  r_{1}\right)  <\xi$. For $k=1,2,\cdots$, we define through
Lemma \ref{Lemma global solution is oscillating},%
\[
r_{2k}=\sup\left\{  r>r_{2k-1}\text{: }h^{\prime}\left(  s\right)  >0\text{
for all }s\in\left(  r_{2k-1},r\right)  \right\}  ,
\]%
\[
r_{2k+1}=\sup\left\{  r>r_{2k}\text{: }h^{\prime}\left(  s\right)  <0\text{
for all }s\in\left(  r_{2k},r\right)  \right\}  .
\]

\begin{lemma}%
\[
\lim_{k\rightarrow\infty}r_{k}=\infty.
\]

\end{lemma}

\begin{proof}
If it is not true, then we have%
\[
\lim_{k\rightarrow\infty}r_{k}=r^{\ast}%
\]
for some $r^{\ast}>0$. Since $h$ is smooth, we have
\[
h\left(  r^{\ast}\right)  =\xi,\quad h^{\prime}\left(  r^{\ast}\right)  =0,
\]
hence Lemma \ref{Lemma global solution is oscillating} implies $h\equiv\xi$,
which is a contradiction.
\end{proof}

Next, we show that the lengths of oscillating intervals $r_{k+1}-r_{k}$ are bounded.

\begin{lemma}
There exists positive constants $C_{1}$ and $C_{2}$ such that%
\[
C_{1}\leq\left(  r_{k+1}-r_{k}\right)  \leq C_{2}%
\]
for any $k=1,2,3,\cdots$.
\end{lemma}

\begin{proof}
Since%
\[
\lim_{k\rightarrow\infty}r_{k}=\infty,
\]
we only need to prove the lemma when $r_{k}$ is sufficiently large.
Differentiating $\left(  \ref{Equation radially symmetric}\right)  $, we have%
\[
h^{\prime\prime\prime}+\frac{N-1}{r}h^{\prime\prime}-\frac{N-1}{r^{2}%
}h^{\prime}+h^{-\alpha-1}h^{\prime}=0.
\]
Let%
\[
w\left(  r\right)  =r^{\frac{N-1}{2}}h^{\prime}\left(  r\right)  ,
\]
then $w$ satisfies%
\begin{equation}
w^{\prime\prime}+\left(  h^{-\alpha-1}-\frac{N^{2}-1}{4r^{2}}\right)  w=0.
\label{Equation derivative of h}%
\end{equation}
Since $h$ is bounded away from both zero and infinity when $r\rightarrow
\infty$ by $\left(  \ref{Equation bounds on h}\right)  $, we have for some
$R>0$ such that for any $r>R$,%
\[
c_{1}\leq h^{-\alpha-1}\left(  r\right)  -\frac{N^{2}-1}{4r^{2}}\leq
c_{2}\text{ }%
\]
for some positive constants $c_{1}$, $c_{2}$. Since $r_{k}$, $k=1,2,\cdots$,
are zeros of $w$, the conclusion follows from Lemma \ref{Lemma Sturm}.
\end{proof}

Finally, we have

\begin{lemma}
\label{Lemma limit of h at infinity}%
\[
\lim_{r\rightarrow\infty}h\left(  r\right)  =\xi.
\]

\end{lemma}

\begin{proof}
Starting with $r_{1}>0$ such that $h^{\prime}\left(  r_{1}\right)  =0$,
$h\left(  r_{1}\right)  <\xi$, we define $r_{k}$ as above. Since $F\left(
h\left(  r_{k}\right)  \right)  $ is monotone decreasing in $k$, and $h\left(
r_{2k}\right)  >\xi$, $h\left(  r_{2k-1}\right)  <\xi$, the property of
function $F$ implies $h\left(  r_{2k}\right)  $ is monotone decreasing and
$h\left(  r_{2k-1}\right)  $ is monotone increasing. Hence we have the limits%
\[
\eta_{1}\equiv\lim_{k\rightarrow\infty}h\left(  r_{2k}\right)  \geq\xi\geq
\eta_{2}\equiv\lim_{k\rightarrow\infty}h\left(  r_{2k-1}\right)  .
\]
Now%
\begin{align*}
&  \eta_{1}-\eta_{2}\leq\left\vert h\left(  r_{k+1}\right)  -h\left(
r_{k}\right)  \right\vert =\int_{r_{k}}^{r_{k+1}}\left\vert h^{\prime}\left(
r\right)  \right\vert dr\\
\leq &  \left(  r_{k+1}-r_{k}\right)  ^{\frac{1}{2}}\left(  \int_{r_{k}%
}^{r_{k+1}}\left\vert h^{\prime}\right\vert ^{2}dr\right)  ^{\frac{1}{2}}%
\leq\sqrt{C_{2}}\left(  \int_{r_{k}}^{r_{k+1}}\left\vert h^{\prime}\right\vert
^{2}dr\right)  ^{\frac{1}{2}},
\end{align*}
which implies%
\begin{equation}
\int_{r_{k}}^{r_{k+1}}\frac{\left\vert h^{\prime}\right\vert ^{2}}{r}%
dr\geq\frac{1}{r_{k+1}}\int_{r_{k}}^{r_{k+1}}\left\vert h^{\prime}\right\vert
^{2}\geq\frac{\left(  \eta_{1}-\eta_{2}\right)  ^{2}}{C_{2}r_{k+1}}\geq
\frac{\left(  \eta_{1}-\eta_{2}\right)  ^{2}}{2C_{2}^{2}}\int_{r_{k}}%
^{r_{k+1}}\frac{1}{r}dr \label{Equation log}%
\end{equation}
when $k$ is sufficiently large. In the last inequality, we used $r_{k+1}%
\leq2r_{k}$ when $k$ is large. From $\left(  \ref{Equation e1}\right)  $, we
have for any $r>r_{1}$,
\[
\int_{r_{1}}^{r}\frac{\left\vert h^{\prime}\right\vert ^{2}}{r}=\frac{1}%
{N-1}\left[  e_{1}\left(  r_{1}\right)  -e_{1}\left(  r\right)  \right]
\leq\frac{e_{1}\left(  r_{1}\right)  }{N-1}.
\]
Therefore $\frac{\left\vert h^{\prime}\right\vert ^{2}}{r}$ is integrable at
$\infty$. Since $\frac{1}{r}$ is not integrable at $\infty$, $\left(
\ref{Equation log}\right)  $ implies that $\eta_{1}=\eta_{2}=\xi$ and%
\[
\lim_{r\rightarrow\infty}h\left(  r\right)  =\xi\text{.}%
\]

\end{proof}

\begin{corollary}
\label{Corollary distance of radius}%
\begin{equation}
\lim_{k\rightarrow\infty}\left(  r_{k+1}-r_{k}\right)  =\pi\left(  \alpha
p\right)  ^{-\frac{1+\alpha}{2\alpha}}. \label{Equation bounds on rk}%
\end{equation}

\end{corollary}

\begin{proof}
In equation $\left(  \ref{Equation derivative of h}\right)  $, we now have%
\[
\lim_{r\rightarrow\infty}\left(  h^{-\alpha-1}-\frac{N^{2}-1}{4r^{2}}\right)
=\xi^{-\alpha-1}=\left(  \alpha p\right)  ^{\frac{1+\alpha}{\alpha}}.
\]
Hence $\left(  \ref{Equation bounds on rk}\right)  $ follows from Lemma
\ref{Lemma Sturm}.
\end{proof}

\section{Nontrivial Smooth Radially Symmetric Solutions}

This section is devoted to the proof of Theorem \ref{Theorem smooth}.

Given $\eta>0$, we consider $\left(  \ref{Equation radially symmetric}\right)
$ with the initial values%
\[
h\left(  0\right)  =\eta>0,h_{r}\left(  0\right)  =0.
\]
The local existence and uniqueness of a smooth solution is standard since $f$
is smooth when $h$ is bounded away from zero. And such solution is actually a
global solution from Corollary \ref{Corollary global}. For any $\eta\neq\xi$,
without loss of generality, we assume $\eta>\xi$. Since%
\[
\left(  r^{N-1}h^{\prime}\right)  ^{\prime}=-r^{N-1}f\left(  h\right)
=-r^{N-1}\left(  -\frac{1}{\alpha}h^{-\alpha}+p\right)  ,
\]
we have $\left(  r^{N-1}h^{\prime}\right)  ^{\prime}<0$ in $\left(
0,\delta\right)  $ for some small $\delta>0$. From $h^{\prime}\left(
0\right)  =0$, we conclude $h^{\prime}\left(  r\right)  <0$ in $\left(
0,\delta\right)  $. Then we can define
\[
r_{1}=\min\left\{  r>0:h^{\prime}\left(  r\right)  =0\right\}  .
\]
The existence of $r_{1}>0$ is guaranteed by Lemma
\ref{Lemma existence of critical points} with $r_{0}=\frac{\delta}{2}$. From
the analysis in the previous section, $h$ will be oscillating around $\xi$,
and all critical points of $h$ can be listed as $r_{1}<r_{2}<r_{3}<\cdots$,
with
\[
C_{1}\leq r_{k+1}-r_{k}\leq C_{2},
\]
and
\[
\lim_{k\rightarrow\infty}\left(  r_{k+1}-r_{k}\right)  =\pi\left(  \alpha
p\right)  ^{-\frac{1+\alpha}{2\alpha}}.
\]
Hence for any $k\geq1$, $h$ is a nontrivial smooth solution of
\[
\left\{
\begin{array}
[c]{lll}%
\triangle h=\frac{1}{\alpha}h^{-\alpha}-p & \text{in} & B_{r_{k}}\left(
0\right)  ,\\
\frac{\partial h}{\partial\nu}=0 & \text{on} & \partial B_{r_{k}}\left(
0\right)  .
\end{array}
\right.
\]
And all nontrivial smooth radial solutions of $\left(
\ref{Equation prescribe volume}\right)  $, when $\Omega$ is a ball, can be
obtained this way. More precisely, let $\Omega=B_{R}\left(  0\right)  $ for a
given $R>0$, then $\left(  \ref{Equation prescribe volume}\right)  $ has a
nontrivial smooth radial solution if and only if $R=r_{k}^{p,\eta}$ for some
$\eta>0,\eta\neq\xi$ and for some $k\geq1$, here we write $r_{k}=r_{k}%
^{p,\eta}$ to recognize its dependence on $p$ and $\eta$.

Now we recall the uniqueness result of M. Del Pino and G. Hernandez
\cite{delpino}.

\begin{proposition}
\label{Proposition Del Pino}Given $\alpha>1$, there exists $d_{0}>0$, such
that%
\begin{equation}
\left\{
\begin{array}
[c]{lll}%
-d\triangle u+u^{-\alpha}=1 & \text{in} & B_{1}\left(  0\right)  ,\\
\frac{\partial u}{\partial\nu}=0 & \text{on} & \partial B_{1}\left(  0\right)
\end{array}
\right.  \label{Equation diffusion}%
\end{equation}
has no nontrivial radial solution whenever $d\geq d_{0}$.
\end{proposition}

It is easy to verify%
\[
\tilde{h}\left(  x\right)  =\left(  \alpha p\right)  ^{\frac{1}{\alpha}%
}h\left(  r_{1}x\right)
\]
satisfies $\left(  \ref{Equation diffusion}\right)  $ with%
\[
d=\frac{\alpha}{\left(  \alpha p\right)  ^{\frac{\alpha+1}{\alpha}}r_{1}^{2}%
}.
\]
Hence Proposition \ref{Proposition Del Pino} implies%
\[
r_{1}>\left(  \frac{\alpha}{\left(  \alpha p\right)  ^{\frac{\alpha+1}{\alpha
}}d_{0}}\right)  ^{\frac{1}{2}}\equiv\left(  \frac{p_{0}}{p}\right)
^{\frac{\alpha+1}{2\alpha}},
\]
here%
\[
p_{0}=\left(  \frac{1}{\alpha^{\frac{1}{\alpha}}d_{0}}\right)  ^{\frac{\alpha
}{\alpha+1}}=\alpha^{-\frac{1}{\alpha+1}}d_{0}^{-\frac{\alpha}{\alpha+1}}.
\]

In general, $r_{k}$ depends on both $p$ and $\eta$. We refer to Corollary
\ref{Corollary scaling} for the scaling of $r_{k}$ when $p,\eta$ changes.

\begin{lemma}
\label{Lemma r1}For any $\eta>0$, $\eta\neq\xi$, we have%
\[
r_{1}\left(  \eta\right)  \geq\sqrt{\frac{2N\alpha\left(  \eta-\xi\right)
}{\xi^{-\alpha}-\eta^{-\alpha}}}.
\]
In particular,%
\[
\lim_{\eta\rightarrow\infty}r_{1}\left(  \eta\right)  =\infty.
\]

\end{lemma}

\begin{proof}
First we assume $\eta>\xi$. From the definition of $r_{1}\left(  \eta\right)
$ and Lemma \ref{Lemma global solution is oscillating}, we have $h^{\prime
}\left(  r_{1}\right)  =0$, $h\left(  r_{1}\right)  <\xi$ and for any
$r\in\left(  0,r_{1}\right)  $, $0<h\left(  r\right)  <\eta$, $h^{\prime
}\left(  r\right)  <0$. Now%
\[
\left(  r^{N-1}h^{\prime}\right)  ^{\prime}=-r^{N-1}\left(  \frac{\xi
^{-\alpha}}{\alpha}-\frac{h^{-\alpha}}{\alpha}\right)  \geq-r^{N-1}\frac
{\xi^{-\alpha}-\eta^{-\alpha}}{\alpha}.
\]
Integrating from $0$ to $r$, we have%
\[
r^{N-1}h^{\prime}\left(  r\right)  \geq-\frac{r^{N}}{N}\frac{\xi^{-\alpha
}-\eta^{-\alpha}}{\alpha},
\]
i.e.%
\[
h^{\prime}\left(  r\right)  \geq-\frac{r}{N\alpha}\left(  \xi^{-\alpha}%
-\eta^{-\alpha}\right)  .
\]
Integrating again from $0$ to $r_{1}$, we have%
\[
h\left(  r_{1}\right)  -h\left(  0\right)  \geq-\frac{r_{1}^{2}}{2N\alpha
}\left(  \xi^{-\alpha}-\eta^{-\alpha}\right)  ,
\]
hence%
\[
r_{1}\left(  \eta\right)  \geq\sqrt{\frac{2N\alpha\left(  h\left(  0\right)
-h\left(  r_{1}\right)  \right)  }{\xi^{-\alpha}-\eta^{-\alpha}}}\geq
\sqrt{\frac{2N\alpha\left(  \eta-\xi\right)  }{\xi^{-\alpha}-\eta^{-\alpha}}%
}.
\]
The bound when $\eta<\xi$ can be proved similarly.
\end{proof}

\section{Rupture Solutions}

In this section, we will consider radial solutions to $\left(
\ref{Equation main}\right)  $ which are not smooth and prove Theorem
\ref{Theorem rupture}. From Corollary \ref{Corollary global}, we need to
consider $h\in C^{0}\left(  \left[  0,\infty\right)  \right)  $ such that
$h\left(  0\right)  =0$ and $h$ satisfies $\left(
\ref{Equation radially symmetric}\right)  $ in $\left(  0,\infty\right)  $.

First, we check the growth rate of $h$ near the origin.

\begin{lemma}
\label{Lemma rupture growth lower bound}Let $h$ be a radially symmetric
rupture solution, then for any $\delta>0$, there exists positive constant
$c_{1}$ such that%
\[
h\left(  r\right)  \geq c_{1}r^{\frac{2}{\alpha+1}}%
\]
holds for any $r\in\left[  0,\delta\right]  $.
\end{lemma}

\begin{proof}
Since $h$ is positive and smooth away from the origin, we only need to prove
the bound for small $\delta$. Let $\delta>0$ be sufficiently small so that%
\[
\frac{1}{\alpha}h^{-\alpha}\left(  r\right)  -p\geq\frac{1}{2\alpha}%
h^{-\alpha}\left(  r\right)
\]
holds for any $r\in\left(  0,\delta\right]  $. Now%
\begin{equation}
\left(  r^{N-1}h^{\prime}\right)  ^{\prime}=r^{N-1}\left(  \frac{h^{-\alpha}%
}{\alpha}-p\right)  \geq\frac{1}{2\alpha}h^{-\alpha}r^{N-1}
\label{Equation 123}%
\end{equation}
implies $r^{N-1}h^{\prime}$ is monotone increasing in $\left(  0,\delta
\right]  $. Since $h\left(  0\right)  =0$ and $h\left(  r\right)  \ $is
positive away from the origin, there exists a sequence $r_{i}\rightarrow0$
such that $h^{\prime}\left(  r_{i}\right)  >0$. Hence, $h^{\prime}\left(
r\right)  >0$ for any $r\in\left(  0,\delta\right]  $. Integrating $\left(
\ref{Equation 123}\right)  $ from $\varepsilon$ to $r$, and using the fact
that $h$ is increasing, we have%
\[
r^{N-1}h^{\prime}\left(  r\right)  -\varepsilon^{N-1}h^{\prime}\left(
\varepsilon\right)  \geq\frac{1}{2N\alpha}h^{-\alpha}\left(  r\right)  \left(
r^{N}-\varepsilon^{N}\right)  .
\]
Letting $\varepsilon\rightarrow0$, we have%
\[
r^{N-1}h^{\prime}\left(  r\right)  \geq\frac{1}{2N\alpha}h^{-\alpha}\left(
r\right)  r^{N}.
\]
Hence for any $r\in\left(  0,\delta\right]  $,%
\[
\frac{d}{dr}h^{\alpha+1}\left(  r\right)  \geq\frac{\alpha+1}{2N\alpha}r.
\]
Integrating from $0$ to $r$, we have
\[
h^{\alpha+1}\left(  r\right)  \geq\frac{\alpha+1}{4N\alpha}r^{2},
\]
i.e. for any $r\in\left(  0,\delta\right]  $,
\[
h\left(  r\right)  \geq\left(  \frac{\alpha+1}{4N\alpha}\right)  ^{\frac
{1}{\alpha+1}}r^{\frac{2}{\alpha+1}}.
\]

\end{proof}

\begin{lemma}
\label{Lemma rupture growth upper bound}Let $h$ be a radially symmetric
rupture solution, then we have for some positive constant $c_{2}$,%
\[
h\left(  r\right)  \leq c_{2}r^{\frac{2}{\alpha+1}}%
\]
for any $r\in\left[  0,\infty\right)  $.
\end{lemma}

\begin{proof}
Since $h$ is uniformly bounded at $\infty$, we only need to prove the
inequality near the origin. First we claim%
\begin{equation}
\lim_{r\rightarrow0^{+}}r^{N-1}h^{\prime}\left(  r\right)  =0.
\label{Equation limit is zero}%
\end{equation}
From%
\[
\left(  r^{N-1}h^{\prime}\right)  ^{\prime}=r^{N-1}\left(  \frac{h^{-\alpha}%
}{\alpha}-p\right)  ,
\]
it follows that $r^{N-1}h^{\prime}$ is monotone increasing near the origin.
Thus, if $\left(  \ref{Equation limit is zero}\right)  $ is false, we would
have%
\[
r^{N-1}h^{\prime}\left(  r\right)  \geq c>0
\]
near the origin, hence%
\[
h^{\prime}\left(  r\right)  \geq cr^{1-N}.
\]
Since $r^{1-N}$ is not integrable near zero, the above inequality contradicts
the fact that $h$ is continuous.

Given $\delta>0$, for any $r\in\left(  0,\delta\right)  $,%
\[
\left(  r^{N-1}h^{\prime}\right)  ^{\prime}=r^{N-1}\left(  \frac{h^{-\alpha}%
}{\alpha}-p\right)  \leq\frac{h^{-\alpha}}{\alpha}r^{N-1}\leq\frac
{c_{1}^{-\alpha}}{\alpha}r^{N-1-\frac{2\alpha}{\alpha+1}}%
\]
by Lemma \ref{Lemma rupture growth lower bound}. Integrating from
$\varepsilon$ to $r$, we obtain%
\[
r^{N-1}h^{\prime}\left(  r\right)  -\varepsilon^{N-1}h^{\prime}\left(
\varepsilon\right)  \leq\frac{1}{\alpha\left(  N-\frac{2\alpha}{\alpha
+1}\right)  }c_{1}^{-\alpha}\left(  r^{N-\frac{2\alpha}{\alpha+1}}%
-\varepsilon^{N-\frac{2\alpha}{\alpha+1}}\right)
\]
Letting $\varepsilon\rightarrow0$, we have%
\[
r^{N-1}h^{\prime}\left(  r\right)  \leq\frac{1}{\alpha\left(  N-\frac{2\alpha
}{\alpha+1}\right)  }c_{1}^{-\alpha}r^{N-\frac{2\alpha}{\alpha+1}},
\]
i.e.,%
\[
h^{\prime}\left(  r\right)  \leq\frac{1}{\alpha\left(  N-\frac{2\alpha}%
{\alpha+1}\right)  }c_{1}^{-\alpha}r^{1-\frac{2\alpha}{\alpha+1}}.
\]
Integrating from $0$ to $r$, we have, for any $r\in\left(  0,\delta\right)
$,
\[
h\left(  r\right)  \leq\frac{\alpha+1}{2\alpha\left(  N-\frac{2\alpha}%
{\alpha+1}\right)  }c_{1}^{-\alpha}r^{\frac{2}{\alpha+1}}.
\]

\end{proof}

Lemmas \ref{Lemma rupture growth lower bound} and
\ref{Lemma rupture growth upper bound} imply that $h\left(  r\right)  $ is of
order $r^{^{\frac{2}{\alpha+1}}}$ near the origin. Now we write%
\begin{equation}
h=c^{\ast}\varphi\left(  r\right)  r^{^{\frac{2}{\alpha+1}}},
\label{Equation change function}%
\end{equation}
where%
\[
c^{\ast}=\left[  \frac{2\alpha}{\alpha+1}\left(  N-2+\frac{2}{\alpha
+1}\right)  \right]  ^{-\frac{1}{\alpha+1}}.
\]
Observe that $h=c^{\ast}r^{\frac{2}{\alpha+1}}$ is a solution of%
\[
\triangle h-\frac{1}{\alpha}h^{-\alpha}=0
\]
in $\left(  0,\infty\right)  $. Direct calculation yields%
\begin{equation}
\varphi^{\prime\prime}+\left(  A+1\right)  \frac{\varphi^{\prime}}{r}%
+\frac{g\left(  \varphi\right)  }{r^{2}}+Cr^{^{-\frac{2}{\alpha+1}}}=0
\label{Equation rupture rescaled}%
\end{equation}
where%
\[
A=N-2+\frac{4}{\alpha+1}>0,\quad C=\frac{p}{c^{\ast}}>0
\]
and%
\[
g\left(  \varphi\right)  =\frac{2}{\alpha+1}\left(  N-2+\frac{2}{\alpha
+1}\right)  \left(  \varphi-\varphi^{-\alpha}\right)  .
\]

Lemmas \ref{Lemma rupture growth lower bound} and
\ref{Lemma rupture growth upper bound} imply that%
\begin{equation}
0<\liminf_{r\rightarrow0^{+}}\varphi\left(  r\right)  \leq\limsup
_{r\rightarrow0^{+}}\varphi\left(  r\right)  <\infty.
\label{Equation rescaled bounds near the origin}%
\end{equation}
On the other hand, let $\varphi$ be a positive solution of $\left(
\ref{Equation rupture rescaled}\right)  $ satisfying $\left(
\ref{Equation rescaled bounds near the origin}\right)  $. Then $h$ defined by
$\left(  \ref{Equation change function}\right)  $ is a rupture solution.

Locally, there exists at least one solution of $\left(
\ref{Equation rupture rescaled}\right)  $ with initial values%
\begin{equation}
\varphi\left(  0\right)  =1,\varphi^{\prime}\left(  0\right)  =0\text{.}
\label{Equation initial condition}%
\end{equation}
To see this, we rewrite the equation as%
\[
\varphi^{\prime\prime}+\left(  A+1\right)  \frac{\varphi^{\prime}}{r}%
+\frac{g^{\prime}\left(  1\right)  \left(  \varphi-1\right)  }{r^{2}}%
+\frac{g\left(  \varphi\right)  -g^{\prime}\left(  1\right)  \left(
\varphi-1\right)  }{r^{2}}+Cr^{^{-\frac{2}{\alpha+1}}}=0
\]
Denoting
\[
\psi=\varphi-1,
\]
we have%
\[
\psi^{\prime\prime}+\left(  A+1\right)  \frac{\psi^{\prime}}{r}+\frac
{g^{\prime}\left(  1\right)  \psi}{r^{2}}+\frac{\tilde{g}\left(  \psi\right)
}{r^{2}}+Cr^{^{-\frac{2}{\alpha+1}}}=0,
\]
where%
\[
\tilde{g}\left(  \psi\right)  =g\left(  \psi+1\right)  -g^{\prime}\left(
1\right)  \psi
\]
satisfies $\tilde{g}\left(  0\right)  =0$, $\tilde{g}^{\prime}\left(
0\right)  =0$. Now let $a_{1},a_{2}$ be two numbers satisfying%
\[
a_{1}+a_{2}=A=N-2+\frac{4}{\alpha+1},\quad a_{1}a_{2}=g^{\prime}\left(
1\right)  =2\left(  N-2+\frac{2}{\alpha+1}\right)  ,
\]
then the real parts of $a_{1},a_{2}$ are both positive and it is easy to
verify that%
\[
\psi^{\prime\prime}+\left(  A+1\right)  \frac{\psi^{\prime}}{r}+\frac
{g^{\prime}\left(  1\right)  \psi}{r^{2}}=r^{-a_{2}-1}\left(  r^{a_{2}%
-a_{1}+1}\left(  r^{a_{1}}\psi\right)  _{r}\right)  _{r}.
\]
Hence, we have%
\begin{align*}
\psi &  =-r^{-a_{1}}\int_{0}^{r}\left\{  s^{a_{1}-a_{2}-1}\int_{0}^{s}%
t^{a_{2}+1}\left(  \frac{\tilde{g}\left(  \psi\left(  t\right)  \right)
}{t^{2}}+Ct^{^{-\frac{2}{\alpha+1}}}\right)  dt\right\}  ds\\
&  =-\frac{C}{\left(  a_{1}+\frac{2\alpha}{\alpha+1}\right)  \left(
a_{2}+\frac{2\alpha}{\alpha+1}\right)  }r^{\frac{2\alpha}{\alpha+1}}%
-r^{-a_{1}}\int_{0}^{r}\left\{  s^{a_{1}-a_{2}-1}\int_{0}^{s}t^{a_{2}-1}%
\tilde{g}\left(  \psi\left(  t\right)  \right)  dt\right\}  ds.
\end{align*}
Let
\[
L\psi=-\frac{C}{\left(  a_{1}+\frac{2\alpha}{\alpha+1}\right)  \left(
a_{2}+\frac{2\alpha}{\alpha+1}\right)  }r^{\frac{2\alpha}{\alpha+1}}%
-r^{-a_{1}}\int_{0}^{r}\left\{  s^{a_{1}-a_{2}-1}\int_{0}^{s}t^{a_{2}-1}%
\tilde{g}\left(  \psi\left(  t\right)  \right)  dt\right\}  ds,
\]
then for $\delta$ sufficiently small, $L$ is a contraction mapping from
\[
X=\left\{  \psi\in C\left(  \left[  0,\delta\right]  \right)  :\left\vert
\psi\left(  r\right)  \right\vert \leq\delta\text{ for any }r\in\left[
0,\delta\right]  \right\}
\]
into itself. Here $L$ is a real mapping even though $a_{1},a_{2}$ could be
complex numbers. Let $\psi$ be the unique fixed point of $L$ in $X$, then
$\varphi=1+\psi$ is a solution to $\left(  \ref{Equation rupture rescaled}%
\right)  $ satisfying $\left(  \ref{Equation initial condition}\right)  $.

Let $\varphi$ be the local solution of $\left(
\ref{Equation rupture rescaled}\right)  $ we just constructed, then $h$
defined by $\left(  \ref{Equation change function}\right)  $ is continuous
with $h\left(  0\right)  =0$ and satisfies $\left(
\ref{Equation radially symmetric}\right)  $ in $\left(  0,\delta\right)  $.
Such solution can be uniquely extended to a solution in $\left(
0,\infty\right)  $ which converges to $\xi$ by Lemma
\ref{Lemma limit of h at infinity}. Thus we have constructed a global rupture solution.

\begin{remark}
\label{Remark weak}From the bounds in Lemmas
\ref{Lemma rupture growth lower bound} and
\ref{Lemma rupture growth upper bound}, it is easy to see that the rupture
solution we constructed is actually a weak solution of $\left(
\ref{Equation main}\right)  $ in $\mathbb{R}^{N}$, $N\geq2$. More precisely,
we have
\[
h\in W_{loc}^{2,P}\left(  \mathbb{R}^{N}\right)  \text{, }h^{-\alpha}\in
L_{loc}^{P}\left(  \mathbb{R}^{N}\right)  \text{ for any }1\leq P<\frac
{\alpha+1}{2\alpha}N,
\]
and
\[
\int_{\mathbb{R}^{N}}h\triangle\phi=\int_{\mathbb{R}^{N}}\left(  \frac
{1}{\alpha}h^{-\alpha}-p\right)  \phi
\]
holds for any $\phi\in C_{0}^{\infty}\left(  \mathbb{R}^{N}\right)  $.
\end{remark}

From the proof of Lemma \ref{Lemma rupture growth lower bound}, we have
$h^{\prime}>0$ near the origin, so we can define%
\[
r_{1}=\min\left\{  r>0\text{: }h^{\prime}\left(  r\right)  =0\right\}  ,
\]
the existence of $r_{1}$ is guaranteed by Lemma
\ref{Lemma existence of critical points}. Furthermore, as in the smooth
solution case, we have a sequence $\left\{  r_{k}\right\}  _{k=1}^{\infty}$
such that for each $k$, $h$ is a rupture solution of%
\[
\left\{
\begin{array}
[c]{lll}%
\triangle h=\frac{1}{\alpha}h^{-\alpha}-p & \text{in} & B_{r_{k}}\left(
0\right)  ,\\
\frac{\partial h}{\partial\nu}=0 & \text{on} & \partial B_{r_{k}}\left(
0\right)  .
\end{array}
\right.
\]

In the remaining part of this section, we will show that the rupture solution
to $\left(  \ref{Equation radially symmetric}\right)  $ is actually unique.

In $\left(  \ref{Equation rupture rescaled}\right)  $, with%
\[
r=e^{-t},\quad\phi\left(  t\right)  =\varphi\left(  r\right)  ,
\]
direct calculation yields%
\begin{equation}
\phi_{tt}-A\phi_{t}+g\left(  \phi\right)  +Ce^{^{-\frac{2\alpha}{\alpha+1}t}%
}=0 \label{equation normalized}%
\end{equation}
on $\left(  -\infty,\infty\right)  $.

\begin{lemma}
\label{Lemma uniqueness of phi}There exists a unique global solution to
$\left(  \ref{equation normalized}\right)  $ satisfying%
\begin{equation}
0<\liminf_{t\rightarrow\infty}\phi\left(  t\right)  \leq\limsup_{t\rightarrow
\infty}\phi\left(  t\right)  <\infty.
\label{Equation normalized rupture bound}%
\end{equation}

\end{lemma}

Noticing that $r\rightarrow0^{+}$ is equivalent to $t\rightarrow\infty$, the
uniqueness of rupture solution follows from Lemma
\ref{Lemma uniqueness of phi}. Before proving Lemma
\ref{Lemma uniqueness of phi}, we first study the behavior of $\phi$ at
$\infty$.

We write
\[
G\left(  \phi\right)  =\frac{2}{\alpha+1}\left(  N-2+\frac{2}{\alpha
+1}\right)  \left(  \frac{\phi^{2}}{2}+\frac{\phi^{1-\alpha}}{\alpha
-1}\right)  ,
\]
hence,%
\[
G^{\prime}\left(  \phi\right)  =g\left(  \phi\right)  .
\]
Multiplying equation $\left(  \ref{equation normalized}\right)  $ with
$\phi_{t}$, we have%
\begin{equation}
\frac{d}{dt}\left(  \frac{\phi_{t}^{2}}{2}+G\left(  \phi\right)  \right)
=A\phi_{t}^{2}-Ce^{^{-\frac{2\alpha}{\alpha+1}t}}\phi_{t}\geq\frac{A}{2}%
\phi_{t}^{2}-\frac{C^{2}}{2A}e^{-\frac{4\alpha}{\alpha+1}t}.
\label{Equation energy inequality}%
\end{equation}

\begin{lemma}
\label{Lemma derivative is bounded}Let $\phi$ be a global solution to $\left(
\ref{equation normalized}\right)  $ satisfying $\left(
\ref{Equation normalized rupture bound}\right)  $, then%
\[
-\infty<\liminf_{t\rightarrow\infty}\phi_{t}\left(  t\right)  \leq
\limsup_{t\rightarrow\infty}\phi_{t}\left(  t\right)  <\infty
\]
and%
\[
\int_{0}^{\infty}\phi_{t}^{2}<\infty.
\]
Furthermore, the limit%
\[
\lim_{t\rightarrow\infty}\left(  \frac{\phi_{t}^{2}}{2}+G\left(  \phi\right)
\right)
\]
exists and is finite.
\end{lemma}

\begin{proof}
If $\phi_{t}$ is unbounded at $\infty$, then $\frac{\phi_{t}^{2}}{2}+G\left(
\phi\right)  $ will be unbounded, so there exists a sequence $\left\{
t_{k}\right\}  _{k=1}^{\infty}$ with $\lim_{k\rightarrow\infty}t_{k}=\infty$,
such that%
\[
\lim_{k\rightarrow\infty}\left(  \frac{\phi_{t}^{2}\left(  t_{k}\right)  }%
{2}+G\left(  \phi\left(  t_{k}\right)  \right)  \right)  =\infty.
\]
For any $t>t_{k}$, integrating%
\[
\frac{d}{dt}\left(  \frac{\phi_{t}^{2}}{2}+G\left(  \phi\right)  \right)
\geq-\frac{C^{2}}{2A}e^{-\frac{4\alpha}{\alpha+1}t},
\]
from $t_{k}$ to $t$, we have%
\begin{align*}
&  \frac{\phi_{t}^{2}\left(  t\right)  }{2}+G\left(  \phi\left(  t\right)
\right)  \geq\frac{\phi_{t}^{2}\left(  t_{k}\right)  }{2}+G\left(  \phi\left(
t_{k}\right)  \right)  -\int_{t_{k}}^{t}\frac{C^{2}}{2A}e^{-\frac{4\alpha
}{\alpha+1}s}ds\\
\geq &  \frac{\phi_{t}^{2}\left(  t_{k}\right)  }{2}+G\left(  \phi\left(
t_{k}\right)  \right)  -\frac{\left(  \alpha+1\right)  C^{2}}{8\alpha A}.
\end{align*}
Hence,%
\[
\lim_{t\rightarrow\infty}\left(  \frac{\phi_{t}^{2}}{2}+G\left(  \phi\right)
\right)  =\infty.
\]
From $\left(  \ref{Equation normalized rupture bound}\right)  $, $G\left(
\phi\right)  $ is bounded at $\infty$, so we deduce%
\[
\lim_{t\rightarrow\infty}\frac{\phi_{t}^{2}}{2}=\infty,
\]
which is impossible for bounded $\phi$. The $L^{2}\left(  0,\infty\right)  $
bound of $\phi_{t}$ follows from $\left(  \ref{Equation energy inequality}%
\right)  $ and the fact that $\frac{\phi_{t}^{2}}{2}+G\left(  \phi\right)  $
is bounded at $\infty$. Finally, since the right hand side of%
\[
\frac{d}{dt}\left(  \frac{\phi_{t}^{2}}{2}+G\left(  \phi\right)  \right)
=A\phi_{t}^{2}-Ce^{^{-\frac{2\alpha}{\alpha+1}t}}\phi_{t}%
\]
is absolutely integrable at $\infty$, we have for any $t_{0}$,%
\[
\lim_{t\rightarrow\infty}\left[  \frac{\phi_{t}^{2}}{2}+G\left(  \phi\right)
\right]  =\frac{\phi_{t}^{2}\left(  t_{0}\right)  }{2}+G\left(  \phi\left(
t_{0}\right)  \right)  +\int_{t_{0}}^{\infty}\left(  A\phi_{t}^{2}%
-Ce^{^{-\frac{2\alpha}{\alpha+1}t}}\phi_{t}\right)
\]
which is finite.
\end{proof}

\begin{lemma}
\label{Lemma one is only possible limit}If%
\[
\lim_{t\rightarrow\infty}\phi=\varsigma,
\]
then
\[
\varsigma=1.
\]

\end{lemma}

\begin{proof}
If $\varsigma\neq1$, we first assume $\varsigma<1$, then for some small
$\delta>0$,
\[
\phi_{tt}-A\phi_{t}=-g\left(  \phi\right)  -Ce^{-\frac{2\alpha}{\alpha+1}%
t}>\delta,
\]
for any $t\geq T_{0}$, where $T_{0}$ is a sufficiently large constant. Hence%
\begin{equation}
\left(  e^{-At}\phi_{t}\right)  _{t}>\delta e^{-At} \label{Equation
small bound}%
\end{equation}
for any $t\geq T_{0}$. Now since $\phi_{t}^{2}$ is integrable, we can choose
$T_{1}>T_{0}$ with $\left\vert \phi_{t}\left(  T_{1}\right)  \right\vert \ $
sufficiently small. For any $t>T_{1}$, integrating $\left(
\ref{Equation small bound}\right)  $ from $T_{1}$ to $t$, we have
\[
e^{-At}\phi_{t}\left(  t\right)  -e^{-AT_{1}}\phi_{t}\left(  T_{1}\right)
>\frac{\delta}{A}\left(  e^{-AT_{1}}-e^{-At}\right)  .
\]
Hence%
\[
\phi_{t}\left(  t\right)  >\left[  \frac{\delta}{A}\left(  e^{-AT_{1}}%
-e^{-At}\right)  +e^{-AT_{1}}\phi_{t}\left(  T_{1}\right)  \right]
e^{At}>\frac{\delta}{2A}e^{-AT_{1}}e^{At}%
\]
when $t$ is sufficiently large, which contradicts the boundedness of $\phi
_{t}$ at $\infty$. The case $\varsigma>1$ can be treated in the same manner.
\end{proof}

\begin{lemma}%
\[
\lim_{t\rightarrow\infty}G\left(  \phi\left(  t\right)  \right)  =G\left(
1\right)  ,
\]
and hence%
\[
\lim_{t\rightarrow\infty}\phi\left(  t\right)  =1.
\]

\end{lemma}

\begin{proof}
Since $G\left(  1\right)  =\min G\left(  \phi\right)  $, if
\[
\lim_{t\rightarrow\infty}\left(  \frac{\phi_{t}^{2}}{2}+G\left(  \phi\right)
\right)  =G\left(  1\right)  ,
\]
then $\lim_{t\rightarrow\infty}G\left(  \phi\right)  $ exists and equals
$G\left(  1\right)  $, as desired. We proceed by contradiction and assume that%
\[
\lim_{t\rightarrow\infty}\left(  \frac{\phi_{t}^{2}}{2}+G\left(  \phi\right)
\right)  =L>G\left(  1\right)  .
\]
We claim%
\[
\liminf_{t\rightarrow\infty}G\left(  \phi\right)  <L.
\]
Otherwise,%
\[
\lim_{t\rightarrow\infty}G\left(  \phi\right)  =L,
\]
which implies%
\[
\lim_{t\rightarrow\infty}\phi=\varsigma
\]
for some $\varsigma$ with $G\left(  \varsigma\right)  =L$, a contradiction to
Lemma \ref{Lemma one is only possible limit}. Hence, there exists a sequence
$\left\{  t_{k}\right\}  _{k=1}^{\infty}$ such that $t_{k}\rightarrow\infty$
and%
\[
G\left(  1\right)  \leq G\left(  \phi\left(  t_{k}\right)  \right)  <L-\delta
\]
for some $\delta>0$. Now we consider
\[
s_{k}=\sup\left\{  s>t_{k}:\text{For any }t\in\left(  t_{k},s\right)  \text{,
}G\left(  \phi\left(  t\right)  \right)  <L-\frac{\delta}{2}\right\}  .
\]
Observe that $s_{k}$ is finite. Otherwise $\phi_{t}^{2}>\frac{\delta}{8}$ for
any $t$ sufficiently large, and then $\phi$ is monotone with derivative
bounded away from zero, hence it will be unbounded, which gives a
contradiction. Since $G\left(  \phi\left(  s_{k}\right)  \right)
=L-\frac{\delta}{2}$, we must have%
\[
\left\vert \phi\left(  t_{k}\right)  -\phi\left(  s_{k}\right)  \right\vert
>\delta_{1},
\]
where $\delta_{1}>0$ is a constant depending on $\delta$, $L$ and $G$. Since
$\phi_{t}$ is bounded, we have
\[
\left\vert s_{k}-t_{k}\right\vert >\delta_{2}\equiv\frac{\delta_{1}%
}{\left\Vert \phi_{t}\right\Vert _{L^{\infty}}}>0.
\]
However, for $t_{k}$ so large that%
\[
\frac{\phi_{t}^{2}}{2}+G\left(  \phi\right)  >L-\frac{\delta}{4}%
\]
on $\left(  t_{k},\infty\right)  $, we have%
\[
\phi_{t}^{2}\left(  t\right)  >\frac{\delta}{2}%
\]
on $\left(  t_{k},s_{k}\right)  $. Hence%
\[
\int_{t_{k}}^{s_{k}}\phi_{t}^{2}>\frac{\delta}{2}\delta_{2}.
\]
On the other hand, Lemma \ref{Lemma derivative is bounded} says%
\[
\lim_{k\rightarrow\infty}\int_{t_{k}}^{\infty}\phi_{t}^{2}=0,
\]
which is a contradiction.
\end{proof}

Now we are ready to prove Lemma \ref{Lemma uniqueness of phi}.

\begin{proof}
[Proof of Lemma \ref{Lemma uniqueness of phi}]Let $\phi$ and $\tilde{\phi}$ be
two global solutions of $\left(  \ref{equation normalized}\right)  $
satisfying $\left(  \ref{Equation normalized rupture bound}\right)  $. Letting
$\psi=\phi-\tilde{\phi}$, we have
\[
\psi_{tt}-A\psi_{t}+B\left(  t\right)  \psi=0.
\]
Here%
\[
A=N-2+\frac{4}{\alpha+1},
\]
and
\[
B\left(  t\right)  =\frac{2}{\alpha+1}\left(  N-2+\frac{2}{\alpha+1}\right)
\left(  1+\frac{\left(  \tilde{\phi}^{-\alpha}-\phi^{-\alpha}\right)  }%
{\phi-\tilde{\phi}}\right)  .
\]
Since%
\[
\lim_{t\rightarrow\infty}\tilde{\phi}\left(  t\right)  =\lim_{t\rightarrow
\infty}\phi\left(  t\right)  =1,
\]
we have%
\[
\lim_{t\rightarrow\infty}B\left(  t\right)  =B_{0}=2\left(  N-2+\frac
{2}{\alpha+1}\right)  >0.
\]
It is easy to check that for any $\lambda$ such that%
\[
\lambda^{2}-A\lambda+B_{0}=0,
\]
we have $\operatorname{Re}\lambda>0$. Since $\psi$ is bounded at $\infty$,
Lemma \ref{Lemma slower decay imply trivial} below with $\lambda_{0}=0$
implies $\psi\equiv0$.
\end{proof}

The following result seems standard and should be well-known. A proof is
included here for the convenience of the reader.

\begin{lemma}
\label{Lemma slower decay imply trivial}Let $u$ satisfy the linear equation%
\begin{equation}
u_{tt}-Au_{t}+B\left(  t\right)  u=0. \label{Equation linear standard}%
\end{equation}
Here $A$ is a constant, and $B\left(  t\right)  $ is a continuous function
such that
\[
\lim_{t\rightarrow\infty}B\left(  t\right)  =B_{0}.
\]
Let $\lambda_{1}$, $\lambda_{2}$ be solutions of%
\[
\lambda^{2}-A\lambda+B_{0}=0.
\]
Suppose that there exists a constant $\lambda_{0}$ satisfying%
\[
\lambda_{0}<\lambda_{m}=\min\left(  \operatorname{Re}\lambda_{1}%
,\operatorname{Re}\lambda_{2}\right)
\]
such that, for some positive constants $T$ and $c$,%
\[
\left\vert u\left(  t\right)  \right\vert \leq ce^{\lambda_{0}t}%
\]
holds for any $t\geq T$, then $u\equiv0$.
\end{lemma}

\begin{proof}
Let $u$ be any function satisfying $\left(  \ref{Equation linear standard}%
\right)  $. For any $\lambda\in\left(  \lambda_{0},\lambda_{m}\right)  $, let
$v=e^{-\lambda t}u$. It is easy to check%
\begin{equation}
v_{tt}-\left(  A-2\lambda\right)  v_{t}+\left(  B\left(  t\right)
-A\lambda+\lambda^{2}\right)  v=0. \label{Equation scaled linear standard}%
\end{equation}
Since $\lambda<\lambda_{m}$, we have%
\[
A-2\lambda=\operatorname{Re}\lambda_{1}+\operatorname{Re}\lambda_{2}%
-2\lambda\geq2\lambda_{m}-2\lambda>0
\]
and%
\[
B_{0}-A\lambda+\lambda^{2}>0.
\]
Multiplying $\left(  \ref{Equation scaled linear standard}\right)  $ with
$v_{t}$, we obtain%
\[
\frac{d}{dt}\left(  v_{t}^{2}+\left(  B_{0}-A\lambda+\lambda^{2}\right)
v^{2}\right)  =2\left(  A-2\lambda\right)  v_{t}^{2}+2\left(  B_{0}-B\left(
t\right)  \right)  vv_{t}.
\]
For any $\varepsilon_{1}>0$ since
\[
\lim_{t\rightarrow\infty}B\left(  t\right)  =B_{0},
\]
there exists $T_{1}>0$, such that
\[
\left\vert \left(  B_{0}-B\left(  t\right)  \right)  vv_{t}\right\vert
\leq\varepsilon_{1}\left(  v_{t}^{2}+\left(  B_{0}-A\lambda+\lambda
^{2}\right)  v^{2}\right)
\]
holds for any $t\geq T_{1}$. Hence for any $t\geq T_{1}$, we have%
\[
\frac{d}{dt}\left(  v_{t}^{2}+\left(  B_{0}-A\lambda+\lambda^{2}\right)
v^{2}\right)  \leq2\left(  A-2\lambda+\varepsilon_{1}\right)  \left(
v_{t}^{2}+\left(  B_{0}-A\lambda+\lambda^{2}\right)  v^{2}\right)  .
\]
Gronwall's inequality then implies that for any $t\geq T_{1}$,
\begin{equation}
v_{t}^{2}+\left(  B_{0}-A\lambda+\lambda^{2}\right)  v^{2}\leq c_{\varepsilon
_{1}}e^{2\left(  A-2\lambda+\varepsilon_{1}\right)  t} \label{Equation v_t}%
\end{equation}
where
\[
c_{\varepsilon_{1}}=\left[  v_{t}^{2}\left(  T_{1}\right)  +\left(
B_{0}-A\lambda+\lambda^{2}\right)  v^{2}\left(  T_{1}\right)  \right]
e^{-2\left(  A-2\lambda+\varepsilon_{1}\right)  T_{1}}.
\]
Now let $u_{1}$ be the solution of $\left(  \ref{Equation linear standard}%
\right)  $ in the Lemma such that%
\[
\left\vert u_{1}\left(  t\right)  \right\vert \leq ce^{\lambda_{0}t}%
\]
holds for any $t\geq T$. Then for any $\lambda\in\left(  \lambda_{0}%
,\lambda_{m}\right)  $ and for any $\varepsilon\in\left(  0,A-2\lambda\right)
$, $v_{1}=e^{-\lambda t}u_{1}$ satisfies%
\begin{align*}
&  \frac{d}{dt}\left(  v_{1,t}^{2}+\left(  B_{0}-A\lambda+\lambda^{2}\right)
v_{1}^{2}\right) \\
=  &  2\left(  A-2\lambda\right)  v_{1,t}^{2}+2\left(  B_{0}-B\left(
t\right)  \right)  v_{1}v_{1,t}\\
\geq &  2\left(  A-2\lambda-\varepsilon\right)  v_{1,t}^{2}-v_{1}^{2}%
\end{align*}
for any $t\geq T_{\varepsilon}$ if we choose $T_{\varepsilon}\geq T$
sufficiently large. Since $A-2\lambda-\varepsilon>0$ and
\begin{equation}
\left\vert v_{1}\left(  t\right)  \right\vert \leq ce^{-\left(  \lambda
-\lambda_{0}\right)  t} \label{Equation v1}%
\end{equation}
holds for any $t\geq T$, we have, by similar arguments as in the proof of
Lemma \ref{Lemma derivative is bounded},
\[
\left\vert v_{1,t}\left(  t\right)  \right\vert \leq C_{1}%
\]
holds for some positive constant $C_{1}$ and for any $t\geq T_{\varepsilon}$.
Hence for any $\lambda\in\left(  \lambda_{0},\lambda_{m}\right)  $,
\[
\left\vert u_{1,t}\left(  t\right)  \right\vert =\left\vert e^{\lambda
t}v_{1,t}\left(  t\right)  +\lambda e^{\lambda t}v_{1}\left(  t\right)
\right\vert \leq C_{2}e^{\lambda t}%
\]
holds for some positive constant $C_{2}$ and for any $t\geq T_{\varepsilon}$.
Especially, for any $\varepsilon_{2}\in\left(  0,\lambda_{m}-\lambda
_{0}\right)  $, we have%
\[
\left\vert u_{1,t}\left(  t\right)  \right\vert \leq C_{2}e^{\left(
\lambda_{0}+\varepsilon_{2}\right)  t}%
\]
holds for any $t\geq T_{\varepsilon}$ where $C_{2}$ is a large constant. Now
for fixed $\lambda\in\left(  \lambda_{0},\lambda_{m}\right)  $, we have%
\begin{equation}
\left\vert v_{1,t}\right\vert =\left\vert e^{-\lambda t}u_{1,t}-\lambda
e^{-\lambda t}u_{1}\right\vert \leq C_{3}e^{\left(  \lambda_{0}-\lambda
+\varepsilon_{2}\right)  t} \label{Equation v1t}%
\end{equation}
for any $t\geq T_{\varepsilon}$ where $C_{3}$ is a large constant. If $u_{1}$
is a nontrivial solution to $\left(  \ref{Equation linear standard}\right)  $,
then $v_{1}$ is a nontrivial solution to $\left(
\ref{Equation scaled linear standard}\right)  $. Let $v_{2}$ be another
solution of $\left(  \ref{Equation scaled linear standard}\right)  $ which is
linearly independent of $v_{1}$. Then $\left(  \ref{Equation v_t}\right)  $
holds for $v_{2}$. Combining with $\left(  \ref{Equation v1}\right)  $ and
$\left(  \ref{Equation v1t}\right)  $, we have for any $t\geq\max\left\{
T_{1},T_{\varepsilon}\right\}  $,%
\[
W\left(  t\right)  =v_{1}v_{2,t}-v_{2}v_{1,t}%
\]
satisfies%
\begin{align}
&  \left\vert W\left(  t\right)  \right\vert \leq\left(  v_{1,t}^{2}+v_{1}%
^{2}\right)  ^{\frac{1}{2}}\left(  v_{2,t}^{2}+v_{2}^{2}\right)  ^{\frac{1}%
{2}}\label{Equation Wt}\\
\leq &  C_{4}e^{\left(  \lambda_{0}-\lambda+\varepsilon_{2}\right)  t}\cdot
e^{\left(  A-2\lambda+\varepsilon_{1}\right)  t}=C_{4}e^{\left(
A-2\lambda+\lambda_{0}-\lambda+\varepsilon_{1}+\varepsilon_{2}\right)
t}.\nonumber
\end{align}
On the other hand, since%
\[
W^{\prime}\left(  t\right)  =\left(  A-2\lambda\right)  W\left(  t\right)  ,
\]
we have%
\begin{equation}
W\left(  t\right)  =W\left(  0\right)  e^{\left(  A-2\lambda\right)  t}
\label{Equation W}%
\end{equation}
Choosing $\varepsilon_{1}$ and $\varepsilon_{2}$ small enough so that
\[
\lambda_{0}-\lambda+\varepsilon_{1}+\varepsilon_{2}<0,
\]
we conclude from $\left(  \ref{Equation Wt}\right)  $ and $\left(
\ref{Equation W}\right)  $ that $W\left(  t\right)  \equiv0$ which contradicts
to the assumption that $v_{1}$, $v_{2}$ are two linearly independent
solutions. Hence $v_{1}\equiv0$ and $u_{1}\equiv0$.
\end{proof}

\section{Scaling of solutions\label{Section scaling}}

In this section, we will use a scaling argument to prove Theorems
\ref{Theorem pk} and \ref{Theorem prescribe volume}.

Let $h^{p,\eta}$ be the unique solution to $\left(
\ref{Equation radially symmetric}\right)  $ satisfying $h\left(  0\right)
=\eta\neq\left(  \alpha p\right)  ^{-\frac{1}{\alpha}}$. When $\eta=0$,
$h^{p,0}$ is the unique rupture solution. Let $r_{k}^{p,\eta}$, $k=1,2,\cdots
$, be the increasing sequence of positive critical points of $h^{p,\eta}$.
Then%
\[
h^{p,\eta,k}\left(  x\right)  =\left(  r_{k}^{p,\eta}\right)  ^{-\frac
{2}{1+\alpha}}h^{p,\eta}\left(  r_{k}^{p,\eta}\left\vert x\right\vert \right)
\]
satisfies%
\[
\left\{
\begin{array}
[c]{lll}%
\triangle h=\frac{1}{\alpha}h^{-\alpha}-p^{p,\eta,k} & \text{in} &
B_{1}\left(  0\right)  ,\\
\frac{\partial h}{\partial\nu}=0 & \text{on} & \partial B_{1}\left(  0\right)
\end{array}
\right.
\]
with%
\[
p^{p,\eta,k}=p\left(  r_{k}^{p,\eta}\right)  ^{\frac{2\alpha}{1+\alpha}}.
\]
Let%
\[
\bar{h}\left(  p,\eta,k\right)  =\frac{1}{\left\vert B_{1}\left(  0\right)
\right\vert }\int_{B_{1}\left(  0\right)  }h^{p,\eta,k}\left(  x\right)
dx=\frac{\left(  r_{k}^{p,\eta}\right)  ^{-\frac{1}{2}}}{\left\vert B_{r_{k}%
}\left(  0\right)  \right\vert }\int_{B_{r_{k}}\left(  0\right)  }h^{p,\eta
}\left(  x\right)  dx.
\]
Then $h^{p,\eta,k}\left(  x\right)  $ is a solution to $\left(
\ref{Equation prescribe volume}\right)  $ with
\[
\bar{h}=\bar{h}\left(  p,\eta,k\right)  .
\]

\begin{lemma}
For any $p>0$ and $\eta\geq0$, $\eta\neq\left(  \alpha p\right)  ^{-\frac
{1}{\alpha}}$,%
\[
h^{p,\eta}\left(  x\right)  =\left(  \alpha p\right)  ^{-\frac{1}{\alpha}%
}h^{\frac{1}{\alpha},\left(  \alpha p\right)  ^{\frac{1}{\alpha}}\eta}\left(
\left(  \alpha p\right)  ^{\frac{1+\alpha}{2\alpha}}x\right)  .
\]

\end{lemma}

\begin{proof}
Let%
\[
f\left(  x\right)  =\left(  \alpha p\right)  ^{-\frac{1}{\alpha}}h^{\frac
{1}{\alpha},\left(  \alpha p\right)  ^{\frac{1}{\alpha}}\eta}\left(  \left(
\alpha p\right)  ^{\frac{1+\alpha}{2\alpha}}x\right)  ,
\]
we have $f\left(  0\right)  =\eta$ and%
\begin{align*}
&  \triangle f\left(  x\right)  =\alpha p\left(  \triangle h^{\frac{1}{\alpha
},\left(  \alpha p\right)  ^{\frac{1}{\alpha}}\eta}\right)  \left(  \left(
\alpha p\right)  ^{\frac{1+\alpha}{2\alpha}}x\right) \\
&  =\alpha p\left(  \frac{1}{\alpha}\left(  h^{\frac{1}{\alpha},\left(  \alpha
p\right)  ^{\frac{1}{\alpha}}\eta}\right)  ^{-\alpha}-\frac{1}{\alpha}\right)
=\frac{1}{\alpha}f^{-\alpha}-p.
\end{align*}
So the lemma follows from the uniqueness of the radial solution.
\end{proof}

\begin{corollary}
\label{Corollary scaling}For each $k$,%
\[
r_{k}^{p,\eta}=\left(  \alpha p\right)  ^{-\frac{1+\alpha}{2\alpha}}%
r_{k}^{\frac{1}{\alpha},\left(  \alpha p\right)  ^{\frac{1}{\alpha}}\eta
},\quad h^{p,\eta,k}=h^{\frac{1}{\alpha},\left(  \alpha p\right)  ^{\frac
{1}{\alpha}}\eta,k},\quad p^{p,\eta,k}=p^{\frac{1}{\alpha},\left(  \alpha
p\right)  ^{\frac{1}{\alpha}}\eta,k}%
\]
and%
\[
\bar{h}\left(  p,\eta,k\right)  =\bar{h}\left(  \frac{1}{\alpha},\left(
\alpha p\right)  ^{\frac{1}{\alpha}}\eta,k\right)  .
\]

\end{corollary}

Now we are ready to prove Theorem \ref{Theorem pk}:

\begin{proof}
[Proof of Theorem \ref{Theorem pk}]Lemma
\ref{Lemma global solution is oscillating} implies that for fixed $\eta>1$,
$r_{k}^{\frac{1}{\alpha},\eta}$, $k=1,2,\cdots$, are well seperated. Hence, we
can apply standard ODE theory to conclude that for each $k=1,2,\cdots$,
$r_{k}^{\frac{1}{\alpha},\eta}$, viewed as a function of $\eta$, is continuous
in $\left(  1,\infty\right)  $. From Lemma \ref{Lemma r1} and noticing
$r_{k}^{\frac{1}{\alpha},\eta}$ is monotone increasing in $k$, we have for
each $k=1,2,\cdots$,
\[
\lim_{\eta\rightarrow\infty}r_{k}^{\frac{1}{\alpha},\eta}=\infty.
\]
Let
\[
R_{k}=\inf_{\eta>1}r_{k}^{\frac{1}{\alpha},\eta},
\]
then $R_{k}$ is monotone nondecreasing in $k$. Furthermore, the interval
$\left(  R_{k},\infty\right)  $ is contained in the range of $r_{k}^{\frac
{1}{\alpha},\eta}$. Given $R>0$, let%
\[
p_{k}=\frac{1}{\alpha}\left(  \frac{R_{k}}{R}\right)  ^{\frac{2\alpha
}{1+\alpha}},
\]
then $p_{k}$ is monotone nondecreasing in $k$. For any $p>p_{k}$, and for any
$1\leq i\leq k$, let $\eta_{i}\in\left(  R_{i},\infty\right)  $ be such that%
\[
p=\frac{1}{\alpha}\left(  \frac{r_{i}^{\frac{1}{\alpha},\eta_{i}}}{R}\right)
^{\frac{2\alpha}{1+\alpha}}.
\]
Then we have%
\[
r_{i}^{p,\left(  \alpha p\right)  ^{-\frac{1}{\alpha}}\eta_{i}}=\left(  \alpha
p\right)  ^{-\frac{1+\alpha}{2\alpha}}r_{i}^{\frac{1}{\alpha},\eta_{i}}=R,
\]
i.e., $h^{p,\left(  \alpha p\right)  ^{-\frac{1}{\alpha}}\eta_{i}}$ is a
nontrivial smooth radial solution to $\left(
\ref{Equation prescribe pressure}\right)  $ in $B_{R}\left(  0\right)  $.
Since $h^{p,\left(  \alpha p\right)  ^{-\frac{1}{\alpha}}\eta_{i}}$, viewed as
a function in $r$, has exactly $i-1$ critical points in $\left(  0,R\right)
$, we have found $k$ distinctive solutions as desired.
\end{proof}

\begin{remark}
For fixed $k\geq1$, numerical computation suggests that $r_{k}^{\frac
{1}{\alpha},\eta}$ is not monotone increasing for $\eta\in\left[  0,1\right)
\cup\left(  1,\infty\right)  $. Hence, given $p>p_{0}$, we may have two
different solutions with the same number of critical points in $\left(
0,R\right)  $.
\end{remark}

From Corollary \ref{Corollary scaling}, to get a solution to $\left(
\ref{Equation prescribe volume}\right)  $ through scaling, we can fix either
$p$ or $\eta$. Without loss of generality, we assume $p=\frac{1}{\alpha}$, and
$\eta\neq1$, $\eta\geq0$.

\begin{proof}
[Proof of Theorem \ref{Theorem prescribe volume}]All radial solutions can be
obtained by scaling. So $\left(  \ref{Equation prescribe volume}\right)  $ has
a rupture solution only when
\[
\bar{h}=\bar{h}_{k}\equiv\bar{h}\left(  \frac{1}{\alpha},0,k\right)  .
\]
Let%
\[
r_{k}\equiv r_{k}^{\frac{1}{\alpha},0},
\]
then%
\[
\bar{h}_{k}=\frac{\left(  r_{k}\right)  ^{-\frac{1}{2}}}{\left\vert B_{r_{k}%
}\left(  0\right)  \right\vert }\int_{B_{r_{k}}\left(  0\right)  }h^{\frac
{1}{\alpha},0}\left(  x\right)  dx.
\]
Since
\[
\lim_{r\rightarrow\infty}h^{\frac{1}{\alpha},0}\left(  r\right)  =1,
\]
we have%
\[
\lim_{k\rightarrow\infty}\frac{1}{\left\vert B_{r_{k}}\left(  0\right)
\right\vert }\int_{B_{r_{k}}\left(  0\right)  }h^{\frac{1}{\alpha},0}\left(
x\right)  dx=1.
\]
Hence the conclusion that
\[
\lim_{k\rightarrow\infty}\sqrt{k\pi}\bar{h}_{k}=1
\]
follows from Corollary \ref{Corollary distance of radius}.
\end{proof}

\section*{Acknowledgements}

We wish to thank Richard Laugesen, Mary Pugh and Dejan Slep\v{c}ev and two
anonymous referees for helpful comments. In particular, we wish to thank Mary
Pugh for pointing out a related work of A. L. Bertozzi, G. Gr\"{u}n and T.P.
Witelski \cite{2001Bertozzi_Grun_Witelski}. The research is supported in part
by the NSF.

\end{document}